\documentclass{amsart}
\usepackage{amsfonts}
\usepackage{amssymb}
 \newtheorem{cor}{Corollary}
 \newtheorem{theorem}{Theorem}
 \newtheorem{lemma}{Lemma}

 \newtheorem{remark}{Remark}
 
 \newcommand{\bee}[1]{\begin{equation}\label{#1}}
 \newcommand{\ene}{\end{equation}}

\begin{document}

\title{ Asymptotics for Capelli Polynomials with Involution}
\author[Benanti]{F. S. Benanti}\address{Dipartimento di Matematica e Informatica
\\ Universit\'a di Palermo \\
via Archirafi, 34\\ 90123
Palermo,Italy}
\email{francescasaviella.benanti@unipa.it}

\author[Valenti]{A. Valenti }\address{Dipartimento di Ingegneria  \\ Universit\`a di Palermo
\\Viale delle Scienze \\
90128 Palermo, Italy}
 \email{angela.valenti@unipa.it}

\thanks{The authors were partially supported by INDAM-GNSAGA of Italy.}
\subjclass[2000]{Primary 16R10; Secondary 16P90} \keywords{Algebras
with involution, Capelli polynomials, Codimension, Growth}

\date{}

\begin{abstract} Let $F\langle X, \ast \rangle$
be the free associative algebra with involution $\ast$ over a field
$F$ of characteristic zero. We study the asymptotic behavior of the
sequence of $\ast$-codimensions of the T-$\ast$-ideal
$\Gamma_{M+1,L+1}^\ast$ of $F\langle X, \ast \rangle$ generated by
the $\ast$-Capelli polynomials $Cap^\ast_{M+1} [Y,X]$ and
$Cap^\ast_{L+1} [Z,X]$ alternanting on $M+1$ symmetric variables and
$L+1$ skew variables, respectively.

 It is well known that, if $F$ is an algebraic closed field of
characteristic zero, every finite dimensional $\ast$-simple algebra
is isomorphic to one of the following algebras:

\begin{itemize}
  \item [$\cdot$]$(M_{k}(F),t)$ the algebra of $k \times k$ matrices with the transpose
involution;
  \item [$\cdot$]$(M_{2m}(F),s)$ the algebra of $2m \times 2m$ matrices with the symplectic
involution;
  \item [$\cdot$]$(M_{h}(F)\oplus M_{h}(F)^{op}, exc)$ the direct sum of the algebra of $h \times h$ matrices and the opposite algebra with
the exchange involution.
\end{itemize}

  We prove that the
$\ast$-codimensions of a finite dimensional $\ast$-simple algebra
are asymptotically equal to the $\ast$-codimensions of
$\Gamma_{M+1,L+1}^\ast$, for some fixed natural numbers $M$ and $L$.
In particular:
$$
c^{\ast}_n(\Gamma^{\ast}_{\frac{k(k+1)}{2} +1,\frac{k(k-1)}{2}
+1})\simeq c^{\ast}_n((M_k(F),t));
$$
$$
c^{\ast}_n(\Gamma^{\ast}_{m(2m-1)+1,m(2m+1)+1})\simeq
c^{\ast}_n((M_{2m}(F),s));
$$
and
$$
c^{\ast}_n(\Gamma^{\ast}_{h^2+1,h^2+1})\simeq
c^{\ast}_n((M_{h}(F)\oplus M_{h}(F)^{op},exc)).
$$

\end{abstract}

\maketitle

\section{Introduction}
Let $(A, \ast)$ be an algebra with involution $\ast$ over a field
$F$ of characteristic zero and let $F\langle X, \ast
\rangle=F\langle x_1,x_1^\ast, x_2, x_2^\ast, \ldots\rangle$ denote
the free associative algebra with involution $\ast$ generated by the
countable set of variables $\{x_1, x_1^\ast, x_2, x_2^\ast,
\ldots\}$ over $F$. Recall that an element $f(x_1,x_1^\ast, \cdots,
x_n,x_n^\ast)$ of $ F\langle X,\ast \rangle$ is a $\ast$-polynomial
identity (or $\ast$-identity) for $A$ if $f(a_1,a_1^\ast, \cdots,
a_n,a_n^\ast)=0$, for all $a_1, \ldots , a_n \in A$. We denote by
$Id^\ast (A)$ the set of all $\ast$-polynomial identities satisfied
by $A$ which is a T-$\ast$-ideal of $F\langle X, \ast \rangle,$
i.e., an ideal invariant under all endomorphisms of $F\langle X,
\ast \rangle$ commuting with the involution of the free algebra. For
$\Gamma=Id^\ast (A)$ we denote by $\mathrm{var}^\ast(\Gamma) =
\mathrm{var}^\ast(A)$ the variety of $\ast$-algebras having the
elements of $\Gamma$ as $\ast$-identities.

It is well known that in characteristic zero $Id^\ast (A)$ is
completely determinated by the multilinear $\ast$-polynomials it
contains. To the T-$\ast$-ideal $\Gamma=Id^\ast (A)$ one can
associates a numerical sequence called the sequence of
$\ast$-codimensions $c_n^\ast(\Gamma)=c_n^\ast(A)$ which is the main
tool for the quantitative investigation of the $\ast$-polynomial
identities of
 $A$. Recall that $c_n^\ast(A)$, $n=1,2,\ldots$,
is the dimension of the space of multilinear polynomial in $n$-th
variables in the corresponding relatively free algebra with
involution of countable rank. Thus, if we denote by  $P_n^\ast$ the
space of all multilinear polynomials of degree $n$ in $x_1,x_1^\ast,
\cdots, x_n,x_n^\ast$ then  $$c_n^\ast(A)= \mathrm{dim} P_n^\ast(A)=
\mathrm{dim }\frac{P_n^\ast}{P_n^\ast\cap Id^\ast (A)}.$$

A celebrated theorem of Amitsur \cite{Am} states that if an algebra
with involution satisfies a $\ast$-polynomial identity then it
satisfies an ordinary  polynomial identity. At the light of this
result in \cite{GR} it was proved that, as in the ordinary case, if
$A$ satisfies a non trivial $\ast$-polynomial identity then
$c_n^\ast(A)$ is exponentially bounded, i.e. there exist constants
$a$ and $b$ such that $c_n^\ast(A)\leq ab^n$, for all $n \geq 1$.
Later (see \cite{BGZ}) an explicit exponential bound for
$c_n^\ast(A)$ was exhibited and in  \cite{GZ0} a characterization of
finite dimensional algebras with involution whose sequence of
$\ast$-codimensions is polynomial bounded was given. This result was
extended to non-finite dimensional algebras (see \cite{GM1}) and
 $\ast$-varieties with almost polynomial growth were
classified in \cite{GM} and \cite{MV}. The asymptotic behavior of
the $\ast$-codimensions was determined in \cite{BGR} in case of
matrices with involution.

Recently (see \cite{GPV}), for any algebra with involution, it was
studied the exponential behavior of $c_n^\ast(A)$, and it was showed
that the $\ast$-exponent of $A$
$$
\textrm{exp}^{\ast}(A)=\displaystyle{\lim_{n\rightarrow \infty}\root n
\of{c_n^{\ast}(A)}}
$$
exists and is a non negative integer. It
should be mentioned that the existence of the $\ast$-exponent was
proved in \cite{GZ00} for finite dimensional algebra with involution.

Now, if $f\in F\langle X, \ast \rangle$ we denote by $\langle f\rangle^\ast$ the T-$\ast$-ideal
generated by $f$. Also for a set of polynomials $V\subset F\langle X, \ast \rangle$ we write $\langle V\rangle^\ast$ to indicate
the T-$\ast$-ideal generated by $V$.

An interesting problem in the theory of PI-algebras with involution
$\ast$  is to describe the T-$\ast$-ideals of $\ast$-polynomial
identities of $\ast$-simple finite dimensional algebras. Recall
that, if $F$ is an algebraically closed field of characteristic
zero, then, up to isomorphisms, all finite dimensional $\ast$-simple
are the following ones  (see \cite{Ro}, \cite{GZ}):

\begin{itemize}
  \item[$\cdot$] $(M_k(F),t)$ the algebra of $ k \times k$
matrices with the transpose involution;
  \item[$\cdot$] $(M_{2m}(F),s)$  the algebra of $ 2m \times 2m
$ matrices with the symplectic involution;
  \item[$\cdot$] $(M_h(F)\oplus M_h(F)^{op},exc)$ the  direct
sum of the algebra of $h \times h$ matrices and the opposite
algebra with the exchange involution.
\end{itemize}

The aim of this paper is to find a relation
among the asymptotics of the $\ast$-codimensions of the finite
dimensional $\ast$-simple algebras and the T-$\ast$-ideals generated
by the $\ast$-Capelli polynomials $Cap^\ast_{M+1} [Y,X]$ and
$Cap^\ast_{L+1} [Z,X]$ alternanting on $M+1$ symmetric variables and
$L+1$ skew variables, respectively.

 More precisely, if $(A, \ast)$ is any algebra with involution $\ast$, let $A^+=\{a
\in A \, | \, a^\ast=a\}$ and $A^-=\{a \in A \, | \, a^\ast=-a\}$
denote the subspaces of symmetric and skew elements of $A$,
 respectively. Since char$F$=0, we can regard the free associative algebra
  with involution $F\langle X, \ast \rangle$ as generated by symmetric and skew variables.
  In particular, for $i=1, 2, \ldots $, we let $y_i=x_i+x_i^\ast$ and $z_i=x_i-x_i^\ast$,
   then we write $X=Y\cup Z$ as the disjoint union of the set $Y$ of symmetric variables
    and the set $Z$ of skew variables and $F\langle X, \ast \rangle=F\langle Y\cup Z\rangle$.
    Hence a polynomial $f=f(y_1, \ldots , y_m,z_1, \ldots , z_n)\in F\langle Y\cup Z\rangle$
     is a $\ast$-polynomial identity of $A$ if and only if $f(a_1, \ldots , a_m,b_1, \ldots , b_n)=0$ for all $a_i \in A^+$, $b_i \in A^-$.

Let us recall that, for any positive integer $m$, the $m$-th Capelli polynomial
 is the element of the free algebra $F\langle X\rangle$ defined as
$$
Cap_m (t_1, \ldots , t_m ; x_1, \ldots , x_{m-1})=
$$
$$
= \sum_{\sigma \in S_m} \mathrm{(sgn
\sigma)}t_{\sigma(1)}x_1t_{\sigma(2)}\cdots t_{\sigma(m-1)}
x_{m-1}t_{\sigma(m)}
$$
where $S_m$ is the symmetric group on $\{1, \ldots , m\}$.
In particular we denote by $$Cap^\ast_m [Y,X]=Cap_m (y_1, \ldots ,
y_m ; x_1, \ldots , x_{m-1})$$
 and $$Cap^\ast_m [Z,X]=Cap_m (z_1, \ldots , z_m ; x_1, \ldots ,
 x_{m-1})$$
 the $m$-th $\ast$-Capelli polynomial in the alternating symmetric
  variables $y_1, \ldots , y_m$ and skew variables $z_1, \ldots , z_m$,
   respectively ($x_1, \ldots ,x_{m-1}$ are arbitrary variables).

Let $Cap_m^+$ denote the set of $2^{m-1}$ polynomials obtained from
$Cap_m [Y,X]$ by deleting any subset of variables $x_i$ (by
evaluating the variables $x_i$ to $1$ in all possible way).
Similarly, we define by  $Cap_m^-$ the set of $2^{m-1}$ polynomials
obtained from $Cap_m [Z,X]$ by deleting any subset of variables
$x_i$.

If $L$ and $M$ are two natural numbers,  we denote  by
$\Gamma_{M+1,L+1}^\ast =\langle  Cap_{M+1}^+, Cap_{L+1}^-  \rangle$ the
T-$\ast$-ideal generated by the polynomials $Cap_{M+1}^+,
Cap_{L+1}^-$. We also write $\mathcal{U}^\ast_{M+1,L+1}$ $=\mathrm{var}^\ast(\Gamma_{M+1,L+1})$ for
the $\ast$-variety generated by $\Gamma_{M+1,L+1}^\ast$.

In this paper we study the asymptotic behavior of the sequence of $
\ast$-codimensions of  $\mathcal{U}^\ast_{M+1,L+1}.$

 Recall that two sequences $a_n$, $b_n$, $n=1,2,\ldots$, are
asymptotically equal, $a_n \simeq b_n$, if $ \lim_{n\rightarrow
+\infty}\frac{a_n}{b_n}=1. $ In the ordinary case (no involution)
(see \cite{GZ1}) it was proved the asymptotic equality between the
codimensions of the Capelli polynomials $Cap_{k^2+1}$ and the
codimensions of the matrix algebra $M_k(F).$ In \cite{Be} these
result was extended to finite dimensional simple superalgebras
proving that the graded codimensions of the $T_2$-ideal generated by
the graded Capelli polynomials $\Gamma_{M+1,L+1},$ for some fixed
$M$, $L$, are asymptotically equal to the graded codimensions of a
simple finite dimensional superalgebra. The link between the
asymptotic of the codimensions of the Amitsur's Capelli-type
polynomials and the verbally prime algebras was studied in
\cite{BS}.

Here we characterize the T-$\ast$-ideal of $\ast$-identities of any
$\ast$-simple finite dimensional algebra showing that

$$
\Gamma^{\ast}_{\frac{k(k+1)}{2} +1,\frac{k(k-1)}{2}
+1}=Id^{\ast}((M_k(F),t)\oplus D');
$$

$$
\Gamma^{\ast}_{m(2m-1)+1,m(2m+1)+1}=Id^{\ast}((M_{2m}(F),s)\oplus
D'');
$$

$$
\Gamma^{\ast}_{h^2+1,h^2+1}=Id^{\ast}((M_{h}(F)\oplus M_{h}(F)^{op},exc)\oplus
D''')
$$

\noindent where $D'$, $D''$ and $D'''$ are  finite dimensional
$\ast$-algebra with exp$^{\ast}(D')<k^2$, exp$^{\ast}(D'')$
$<(2m)^2$ and exp$^{\ast}(D''')<2h^2$. It follows that
asymptotically

$$
c^{\ast}_n(\Gamma^{\ast}_{\frac{k(k+1)}{2} +1,\frac{k(k-1)}{2}
+1})\simeq c^{\ast}_n((M_k(F),t));
$$

$$
c^{^\ast}_n(\Gamma^{\ast}_{m(2m-1)+1,m(2m+1)+1})\simeq
c^{\ast}_n((M_{2m}(F),s));
$$

$$
c^{^\ast}_n(\Gamma^{\ast}_{h^2+1,h^2+1})\simeq
c^{\ast}_n((M_{h}(F)\oplus M_{h}(F)^{op},exc)).
$$

\section{Preliminaries}

 Let $F$ be a field of characteristic zero and let  $G$
be the Grassmann algebra over $F$ generated by the elements $e_1,
e_2, \ldots$ subject to the following condition $e_ie_j = -e_je_i,$
for all $i, j \ge 1$. Recall that $G$ has a natural $Z_2$-grading
$G= G_0 \oplus G_1$ where $G_0$ (resp. $G_1$) is the span of the
monomials in the $e_i's$ of even length (resp. odd length). If $B=
B_0 \oplus B_1$ is  a superalgebra, then the Grassmann envelope of
$B$ is defined as $$G(B) = ( G_0 \otimes B_0)\oplus(G_1 \otimes
B_1).$$ The relevance of $G(A)$ relies in a result of Kemer
(\cite[Theorem 2.3]{Ke}) stating that if $B$ is any PI-algebra, then
its T-ideal of polynomial identities coincides with the T-ideal of
identities of the Grassmann envelope of  a suitable finite
dimensional superalgebra. This result has been extended to algebras
with involution in \cite{alj-gia-kar} and the following result holds

\begin{theorem}
If $A$ is a $PI$-algebra with involution over a field $F$ of characteristic zero, then there exists a finite dimensional superalgebra  with
superinvolution $B$ such that $Id^*(A) = Id^*(G(B))$.
 \end{theorem}

Recall that a superinvolution $*$ of $B$ is  a linear map of $B$ of
order two such that $(ab)^*= (-1)^{\mid a\mid\mid b\mid} b^*a^*,$
for any homogeneous elements $a,b \in B,$ where $\mid a\mid$ denotes
the homogeneous degree of $a$. It is well known that in this case
$B^*_0\subseteq B_0, B^*_1\subseteq B_1$ and we decompose $B=
B_0^+\bigoplus B_0^- \bigoplus B_1^+ \bigoplus B_1^-.$

We can define a superinvolution $*$ on $G$ by requiring that $e_i^*
= -e_i,$ for any $i \ge 1.$ Then it is easily checked that $G_0 =
G^+$ and $ G_1 = G^-.$ Now, if $B$ is  a superalgebra one can
perform its Grassmann envelope $G(B)$ and in  \cite{alj-gia-kar} it
was shown that if $B$ has a superinvolution $*$ we can regard $G(B)$
as an algebra with involution by setting $(g \otimes a)^* = g^*
\otimes a^*,$ for homogeneous elements $g \in G, a \in B.$

 By making use of the previous theorem, in \cite{GPV} it was proved the existence of
 the $\ast$-exponent  of a $PI$-algebra with involution $A$ and also an explicit way of computing
 exp$^\ast(A)$  was given.
More precisely if  $B$  is a finite dimensional algebra with
superinvolution over an algebraic closed field of characteristic
zero, then by \cite{GIL} we write $B= \bar{B}+J$ where $\bar{B}$ is
a maximal semisimple superalgebra with induced superinvolution and
$J =J(B)= J^*.$ Also we can write
$$\bar{B} = B_1 \oplus \cdots \oplus B_k$$
where $B_1, \cdots, B_k$ are simple superalgebras with induced
superinvolution. We say that a subalgebra
  $B_{i_1} \oplus \cdots \oplus B_{i_t},$ where $B_{i_1},\ldots,  B_{i_t}$
  are distinct simple components, is admissible if for some permutation $(l_1, \ldots, l_t)$ of
  $(i_1, \ldots, i_t)$ we have that $B_{l_1}JB_{l_2}J \cdots  JB_{l_t}  \ne 0.$
Moreover if $B_{i_1} \oplus \cdots \oplus B_{i_t}$ is an admissible
subalgebra of
   $B$ then $B' = B_{i_1} \oplus \cdots \oplus B_{i_t}+J$  is
  called a reduced algebra.
 In \cite{GPV} it was proved
that  exp$^\ast(A)$ = exp$^\ast(G(B))= d$ where $d$ is the maximal
dimension of an admissible subalgebra of $B.$

It follows immediately that

 \begin{remark} If $A$ is a $*$-simple algebra then
$exp^\ast(A) = \textrm{dim}_F A$.
\end{remark}

We next prove that the reduced algebras are basic elements of any
$*$-variety. We start with the following

\begin{lemma} \label{codimensioni}
Let $A$ and $B$ be algebras with involution satisfying a
$\ast$-polynomial identity. Then

$$
c_n^{\ast}(A), c_n^{\ast}(B) \leq c_n^{\ast}(A\oplus B) \leq
c_n^{\ast}(A)+c_n^{\ast}(B).
$$

\noindent Hence the $\mathrm{exp}^{\ast}(A\oplus
B)=\mathrm{max}\{\mathrm{exp}^{\ast}(A), \mathrm{exp}^{\ast}(B)\}.$

\end{lemma}

\noindent {\bf Proof.} It follows easily from the proof of the Lemma
1 in \cite{GZ1}.

\bigskip

If $ \mathcal{V}= \mathrm{var}^{\ast}(A)$  is the variety of $\ast$-algebras
generated by $A$ we write  $Id^\ast (\mathcal{V})= Id^{\ast}(A)$,
$c_n^{\ast}(\mathcal{V})=c_n^{\ast}(A)$ and
$\textrm{exp}^{\ast}(\mathcal{V}) = \textrm{exp}^{\ast}(A)$.

We have the following

\begin{theorem} \label{decomposizione}
Let $\mathcal{V}$ be a proper variety of $\ast$-algebras. Then there
exists a finite number of reduced superalgebras with superinvolution
$B_1, \ldots, B_t$ and a finite dimensional superalgebra with
superinvolution $D$ such that

$$
\mathcal{V}=\mathrm{var}(G(B_1)\oplus \cdots \oplus G(B_t) \oplus G(D))
$$

\noindent with
$\mathrm{exp}^{\ast}(\mathcal{V})=\mathrm{exp}^{\ast}(G(B_1))=\cdots=\mathrm{exp}^{\ast}G((B_t))$
and $\mathrm{exp}^{\ast}(G(D))<\mathrm{exp}^{\ast}(\mathcal{V}).$
\end{theorem}

\noindent {\bf Proof.}  The proof follows closely the proofs given
in \cite[Theorem 1]{GZ1} and in \cite[Theorem 3]{Be}.  Let $A$ be a
$\ast$-PI-algebra  such that $\mathcal{V}=\mathrm{var}^\ast(A)$. By
Theorem 1, there exists a finite dimensional superalgebra with
superinvolution $B$ such that $Id^{\ast}(A) = Id^{\ast}(G(B))$.
Also, by \cite[Theorem 4.1]{GIL}, we may assume that

$$
B=\bar{B}_1\oplus\cdots\oplus \bar{B}_s+J(B),
$$
\noindent where  $\bar{B}_i$ are simple $\ast$-superalgebras and
$J^\ast=J$ is the Jacobson radical of $B$. Let
$\textrm{exp}^{\ast}(A)=d$.  Then, since $d$ is the maximal
dimension of an admissible subalgebra of $A$  there exist distinct
$\ast$-simple superalgebras $\bar{B}_{j_1}, \ldots \bar{B}_{j_k} $
such that

$$
\bar{B}_{j_1}J \cdots J\bar{B}_{j_k}\neq 0 \,\,\,\,\,\,\,
\textrm{and} \,\,\,\,\,\,\, \textrm{dim}_F(\bar{B}_{j_1}\oplus
\cdots \oplus \bar{B}_{j_k})=d.
$$

\noindent Let $\Gamma_1, \ldots , \Gamma_t$ be all possible subset
of $\{1, \ldots , s\}$ such that, if  $\Gamma_j=\{j_1, \ldots
,j_k\}$, then $dim_F(\bar{B}_{j_1}\oplus \cdots \oplus
\bar{B}_{j_k})=d$ and $\bar{B}_{\sigma(j_1)}J \cdots
J\bar{B}_{\sigma(j_k)}\neq 0$ for some permutation $\sigma \in S_k$.
For any such $\Gamma_j$, $j=1, \ldots t,$  then we put
$B_j=\bar{B}_{j_1}\oplus \cdots \oplus \bar{B}_{j_k} + J$. It
follows, by the characterization of the $\ast$-exponent, that
$$
\mathrm{exp}^{\ast}(G(B_1))=\cdots
=\mathrm{exp}^{\ast}(G(B_t))=d=\mathrm{exp}^{\ast}(G(B)).
$$
\noindent Let $D=D_1 \oplus \cdots \oplus D_p$, where $D_1, \ldots ,
D_p$ are all subsuperalgebras of $B$ with superinvolution of the
type $\bar{B}_{i_1}\oplus \cdots \oplus \bar{B}_{i_r} + J$, with
$1\leq i_1 <\cdots < i_r \leq s$ and
$\textrm{dim}_F(\bar{B}_{i_1}\oplus \cdots \oplus \bar{B}_{i_r})<d
$. Then, by Lemma \ref{codimensioni}, we have
$\mathrm{exp}^{\ast}(G(D))< \mathrm{exp}^{\ast}(G(B)).$

 Now, we want to prove
that $\mathrm{var}^\ast(G(B_1) \oplus \cdots \oplus G(B_t)\oplus
G(D))=\mathrm{var}^\ast(G(B))$.  Since $G(D), G(B_i)\in
\mathrm{var}^\ast(A)$, $\forall i=1, \ldots , t,$ it follows that

$$
\mathrm{var}^\ast(G(B_1) \oplus \cdots \oplus G(B_t)\oplus
G(D))\subseteq \mathrm{var}^\ast(G(B)).
$$

Now, let $f=f(y_1^+, \ldots,y_n^+,y_1^-, \ldots,y_m^-,z_1^+, \ldots,
z_p^+,z_1^-, \ldots, z_q^-)$ be a multilinear polynomial such that
$f\not\in Id^{\ast}(G(B))$. We shall prove that $f \not\in
Id^{\ast}(G(B_1) \oplus \cdots \oplus G(B_t)\oplus G(D))$. Since
$f\not\in Id^{\ast}(G(B))$, there exist $$a_{1,0}^+\otimes g_{1,0},
\ldots , a_{n,0}^+\otimes g_{n,0} \in G(B)_0^+=B_0^+\otimes G_0,$$
$$a_{1,0}^-\otimes h_{1,0}, \ldots , a_{m,0}^-\otimes h_{m,0} \in
G(B)_0^-=B_0^-\otimes G_0,$$ $$b_{1,1}^-\otimes g_{1,1}, \ldots ,
b_{p,1}^-\otimes g_{p,1} \in G(B)_1^+=B_1^-\otimes G_1,$$ and
$$b_{1,1}^+\otimes h_{1,1}, \ldots , b_{q,1}^+\otimes h_{q,1} \in
G(B)_1^-=B_1^+\otimes G_1$$ such that
$$
f(a_{1,0}^+\otimes g_{1,0}, \ldots , a_{n,0}^+\otimes g_{n,0},
a_{1,0}^-\otimes h_{1,0}, \ldots , a_{m,0}^-\otimes h_{m,0},
$$
$$b_{1,1}^-\otimes g_{1,1}, \ldots , b_{p,1}^-\otimes
g_{p,1},b_{1,1}^+\otimes h_{1,1}, \ldots , b_{q,1}^+\otimes
h_{q,1})\neq 0.
$$
It follows that
$$
0 \neq f(a_{1,0}^+\otimes g_{1,0}, \ldots , a_{n,0}^+\otimes g_{n,0},
a_{1,0}^-\otimes h_{1,0}, \ldots , a_{m,0}^-\otimes h_{m,0},$$
$$
b_{1,1}^-\otimes g_{1,1}, \ldots , b_{p,1}^-\otimes
g_{p,1},b_{1,1}^+\otimes h_{1,1}, \ldots , b_{q,1}^+\otimes
h_{q,1})=
$$
$$
\tilde{f}(a_{1,0}^+,\ldots , a_{n,0}^+,a_{1,0}^-, \ldots ,
a_{m,0}^-,b_{1,1}^-, \ldots ,b_{p,1}^-,b_{1,1}^+, \ldots ,
b_{q,1}^+)\otimes$$
$$ g_{1,0} \cdots g_{n,0}h_{1,0} \cdots
h_{m,0}g_{1,1} \cdots g_{p,1}h_{1,1} \cdots h_{q,1}
$$

\noindent where $\tilde{f}$ is the multilinear polynomial introduced
in \cite[Lemma 1]{GPV}. Clearly $\tilde{f}\neq 0$. From the
linearity of $\tilde{f}$ we can assume that $a_{i,0}^+$,
$b_{j,1}^-$, $a_{i,0}^-$ and $b_{j,1}^+$ $\in B_1 \cup \cdots \cup
B_s \cup J$. Since $B_iB_j=0$ for $i \neq j$, from the property of
the $\ast$-exponent described above, we have that
$$
a_{1,0}^+, \ldots , a_{n,0}^+,a_{1,0}^-, \ldots ,
a_{m,0}^-,b_{1,1}^-, \ldots ,b_{p,1}^-,b_{1,1}^+, \ldots , b_{q,1}^+
\in B_{j_1}\oplus \cdots \oplus B_{j_k} + J
$$
for some $B_{j_1}, \ldots , B_{j_k}$ such that  dim$_F(B_{j_1}\oplus
\cdots \oplus B_{j_k})\leq d.$

 Thus $f$ is not an identity for one
of the algebras $G(B_1), \ldots , G(B_t), G(D)$. Hence $f \not\in
Id^{\ast}(G(B_1) \oplus \cdots \oplus G(B_t)\oplus G(D))$. In
conclusion
$$
\mathrm{var}^\ast(G(B))\subseteq \mathrm{var}^\ast (G(B_1)\oplus
\cdots \oplus G(B_t) \oplus G(D))
$$
and the proof is complete.

\medskip

An application of Theorem \ref{decomposizione} is given in
 terms of $\ast$-codimensions.

\bigskip

\begin{cor} \label{starcodimensioni}
Let $\mathcal{V}=\mathrm{var}^\ast(A)$ be a proper variety of
$\ast$-algebras. Then there exists a finite number of reduced
superalgebras with superinvolution $B_1, \ldots, B_t$ and a finite
dimensional superalgebra with superinvolution $D$ such that
$$
c_n^{\ast}(A)\simeq c_n^{\ast} (G(B_1)\oplus \cdots \oplus G(B_t)).
$$
\end{cor}
\noindent {\bf Proof.} By Theorem \ref{decomposizione}, there is a
finite number of reduced superalgebras with superinvolution $B_1,
\ldots, B_t$ such that
$$
\mathcal{V}=\mathrm{var}^\ast(A)= \mathrm{var}^\ast(G(B_1)\oplus
\cdots \oplus G(B_t)\oplus G(D))
$$

\noindent with
$\textrm{exp}^{\ast}(A)=\textrm{exp}^{\ast}(G(B_1))=\cdots=\textrm{exp}^{\ast}(G(B_t))$
and $\textrm{exp}^{\ast}(G(D))< \textrm{exp}^{\ast}(A)$. Then, by
Lemma \ref{codimensioni}
$$
c_n^{\ast}(G(B_1)\oplus \cdots \oplus G(B_t)) \leq
c_n^{\ast}(G(B_1)\oplus \cdots \oplus G(B_t)\oplus G(D))
\leq
$$
$$
c_n^{\ast}(G(B_1)\oplus \cdots \oplus G(B_t))+c_n^{\ast}(G(D)).
$$

Recalling that
$\textrm{exp}^{\ast}(G(D))<\textrm{exp}^{\ast}(G(B_1))=\mathrm{exp}^{\ast}(G(B_1)\oplus
\cdots \oplus G(B_t))$ we have that

$$
c_n^{\ast}(A)\simeq c_n^{\ast} (G(B_1)\oplus \cdots \oplus G(B_t))
$$

\noindent and the proof of the corollary is complete.
\medskip

\bigskip

If $A$ is a finite dimensional $\ast$-algebra we obtain a simplified
form of the previous theorem and corollary. Let us recall that an
algebra with involution can be regarded as a superalgebra with
superinvolution with trivial grading. We have the following

\bigskip

\begin{cor} \label{decomposizionefinita}
Let $A$ be a finite dimensional $\ast$-algebra. Then there exists a
finite number of reduced $\ast$-algebras  $B_1, \ldots, B_t$ and a
finite dimensional $\ast$-algebra  $D$ such that
$$\mathrm{var}^\ast(A)=\mathrm{var}^\ast(B_1\oplus \cdots \oplus
B_t\oplus D)$$
$$
c_n^{\ast}(A)\simeq c_n^{\ast} (B_1\oplus \cdots \oplus B_t)
$$ and $$\mathrm{exp}^\ast(A)= \mathrm{exp}^\ast(B_1)=\cdots =
\mathrm{exp}^\ast(B_t),\,\,  \mathrm{exp}^\ast(D)<\mathrm{exp}^\ast(A).$$

\end{cor}

\bigskip

The following results give us a characterization of the
$\ast$-varieties satisfying a Capelli identity. Let's start with the
\medskip

\begin{remark} \label{capelliidentity}
Let $M$ and $L$ be two natural numbers. If $A$ is an algebra with
involution satisfying the $\ast$-Capelli polynomials $Cap^\ast_M
[Y,X]$ and $Cap^\ast_L [Z,X]$, then $A$ satisfies the Capelli
identity $Cap_{M+L} (x_1, \ldots , x_{M+L} ; \bar{x}_1, \ldots ,
\bar{x}_{M+L-1})$.
\end{remark}
\noindent {\bf Proof.} To obtain the thesis it is sufficient to
observe that
$$
Cap_{M+L} (x_1, \ldots , x_{M+L} ; \bar{x}_1, \ldots ,
\bar{x}_{M+L-1})=
$$
$$
Cap_{M+L} (\frac{x_{1}+x_{1}^\ast}{2}+\frac{x_{1}-x_{1}^\ast}{2},
\ldots ,
\frac{x_{M+L}+x_{M+L}^\ast}{2}+\frac{x_{M+L}-x_{M+L}^\ast}{2} ;
\bar{x}_1, \ldots , \bar{x}_{M+L-1})
$$

\noindent is a linear combinations of $\ast$-Capelli polynomials
alternating either in $m\geq M$ symmetric variables or in $l\geq L$
skew variables.

\bigskip

The proof of the next result follows closely the proof given in
\cite[Theorem 11.4.3]{GZ}

\bigskip

\begin{theorem} \label{finitelygenerated}
Let $\mathcal{V}$ be a variety of $\ast$-algebras. If $\mathcal{V}$
satisfies the Capelli identity of some rank then
$\mathcal{V}=\mathrm{var}^\ast(A)$, for some finitely generated
$\ast$-algebra $A$.
\end{theorem}

\bigskip
Let $M$, $L$ be two natural numbers. Let $A=A^{+}\oplus A^{-}$ be a
generating $\ast$-algebra of $\mathcal{U}^{\ast}_{M+1,L+1}$. By
remark \ref{capelliidentity}, $A$ satisfies a Capelli identity.
Hence by the previous theorem, we may assume that $A$ is a finitely
generated $\ast$-algebra. Moreover  by \cite[Theorem 1]{Sv}  we may
consider $A$ as a finite-dimensional $\ast$-algebra. Since any
polynomial alternating on $M+1$ symmetric variables vanishes in $A$
(see \cite[Proposition 1.5.5]{GZ}), we get that $\mathrm{dim}\,
A^{+}\leq M$.  Similarly we get that $\mathrm{dim}\, A^{-}\leq L$
and $\mathrm{exp}^\ast(A)\leq \mathrm{dim}\, A \leq M+L$. Thus we
have the following

\bigskip

\begin{lemma}\label{esponente}
$\mathrm{exp}^\ast(\mathcal{U}^{\ast}_{M+1,L+1})\leq M+L$.
\end{lemma}

\bigskip

\section{The algebra $UT^\ast(A_1, \ldots, A_n)$ }

In this section we recall the construction of the $\ast$-algebra
$UT^\ast(A_1, \ldots, A_n)$ given in Section 2 of \cite{DvLs}.
 Let
$A_1, \ldots, A_n$ be a $n$-tuple of finite dimensional $\ast$-simple algebras, then
$A_i=(M_{d_i},\mu_i)$, where $\mu_i$ is the transpose or the
 symplectic involution, or $A_i=(M_{d_i}\oplus M_{d_i}^{op},exc)$,
where $exc$ is the exchange involution.

Let $\gamma_d$ be the orthogonal involution defined on the matrix
algebra $M_d(F)$ by putting, for all $a \in M_d(F)$,
$$a^{\gamma_d}=g^{-1}a^tg=ga^tg,$$ where
$$
g=\left(
    \begin{array}{ccccc}
      0 & & \ldots & & 1 \\
       & & &\cdot & \\
       & &\cdot  & & \\
         &\cdot & & & \\
      1 & &\ldots & & 0 \\
    \end{array}
  \right)
$$
and $a^t$ is the transposed of the matrix $a$. $\gamma_d$ acts on
matrix units $e_{pq}$ of $M_d$ by sending it to $e^{\gamma_d}_{pq} =
e_{d-q+1,d-p+1}$ (it is the reflection along the secondary
diagonal).

\noindent If $d=\sum_{i=1}^n \textrm{dim}_FA_i$, then we have an
embedding of $\ast$-algebras

 $$\Delta : \bigoplus_{i=1}^n A_i \rightarrow
(M_{2d}(F), \gamma_{2d})$$ defined by

$$
(a_1, \ldots, a_n)\rightarrow \left(
                                \begin{array}{cccccc}
                                  \bar{a}_1 &  &  &  &  & \\
                                  & \ddots &  &  &  &  \\
                                  &  &  \bar{a}_n &  &  &  \\
                                  &  &  & \bar{b}_n &  &  \\
                                  &  &  &  & \ddots &  \\
                                  &  &  &  &  & \bar{b}_1 \\
                                \end{array}
                              \right)
$$

\noindent where, if $a_i \in A_i=(M_{d_i},\mu_i)$, then
$\bar{a}_i=a_i$ and $\bar{b}_i=a_i^{\mu_i\gamma_{d_i}}$, and if
$a_i=(\tilde{a}_i, \tilde{b}_i)\in A_i=(M_{d_i}\oplus
M_{d_i}^{op},exc)$, then $\bar{a}_i=\tilde{a}_i$ and
$\bar{b}_i=\tilde{b}_i$.

Let denote by $D=D(A_1, \ldots, A_n)\subseteq M_{2d}(F)$  the
$\ast$-algebra  image of $\bigoplus_{i=1}^n A_i$ by $\Delta$ and let
 $U$ be the subspace of $M_{2d}(F)$ so defined:

$$
\left(
  \begin{array}{cccccccc}
    0 & U_{12} & \cdots & U_{1t} &  &  &  &  \\
      & \ddots & \ddots & \vdots &  &  &  &  \\
      &   & 0 & U_{t-1t} &  &  &  &  \\
      &   &   & 0 &  &  &   &   \\
      &   &   &   & 0 & U_{tt-1} & \cdots & U_{t1} \\
      &   &   &   &   & \ddots & \ddots & \vdots \\
      &   &   &   &   &   & 0 & U_{21} \\
      &   &   &   &   &   &   & 0 \\
  \end{array}
\right)
$$

\noindent where, for $1 \leq i,j \leq n$, $i\neq j$, $U_{ij}$ denote
the vector space of the rectangular matrices of dimensions $d_i
\times d_j.$ Let define

$$
UT^\ast (A_1, \ldots, A_n)=D\oplus U \subseteq M_{2d}(F)
$$

\noindent (see section 2 of \cite{DvLs}).

It is easy to show that $UT^\ast (A_1, \ldots, A_n)$ is a subalgebra
with involution of $(M_{2d}(F),\gamma_{2d})$ in which the algebras
$A_i$ are embedded as $\ast$-algebras and whose $\ast$-exponent is
given by
$$\mathrm{exp}^\ast(UT^\ast (A_1, \ldots, A_n))=\sum_{i=1}^n \mathrm{dim}_FA_i.$$

In \cite{DvS1} and \cite{DvS2} the link between the degrees of
$\ast$-Capelli polynomials and the $\ast$-polynomial identities of
$UT^\ast (A_1, \ldots, A_n)$ was investigated.

If we set $d^+:=\sum_{i=1}^n \mathrm{dim}_FA_i^+$ and $d^-:=\sum_{i=1}^n
\mathrm{dim}_FA_i^-, $ then the following result applies (see \cite{DvS})

\begin{lemma}\label{capelliidentita} Let $R=UT^\ast (A_1, \ldots, A_n).$
Then $Cap^\ast_M [Y,X]$ and $Cap^\ast_L [Z,X]$ are in $Id^\ast(R)$
if and only if $M\geq d^++n$ and $L\geq d^-+n$.
\end{lemma}

\bigskip

\section{Asymptotics for $\mathcal{U}^\ast_{\frac{k(k+1)}{2} +1,\frac{k(k-1)}{2} +1}$ and $(M_k(F),t)$}

Let $A= \bar{A} \oplus J$ where $\bar{A}$ is a  $\ast$-simple finite
dimensional algebra and $J=J(A)$ is its Jacobson radical. It
is well known that the Jacobson radical $J$ is a $\ast$-ideal of
$A$.

We start with the
following key lemmas that hold for any $\ast$-simple finite
dimensional algebra.

\begin{lemma} \label{radicaltrasposta} Let $A= \bar{A} \oplus J$
where $\bar{A}$ is a  $\ast$-simple finite dimensional algebra and
$J=J(A)$ is its Jacobson radical. Then $J$ can be decomposed
into the direct sum of four $\bar{A}$-bimodules
$$ J=J_{00}\oplus J_{01} \oplus J_{10} \oplus J_{11} $$ where, for
$p,q \in \{0,1\}$, $J_{pq}$ is a left faithful module or a $0$-left
module according to $p=1$, or $p=0$, respectively. Similarly,
$J_{pq}$ is a right faithful module or a $0$-right module according
to $q=1$ or $q=0$, respectively. Moreover, for $p,q,i,l \in
\{0,1\}$, $J_{pq}J_{ql}\subseteq J_{pl}$, $J_{pq}J_{il}=0$ for
$q\neq i$ and there exists a finite dimensional nilpotent
$\ast$-algebra $N$ such that $J_{11}\cong \bar{A}\otimes_F N$ (isomorphism
of $\bar{A}$-bimodules and of $\ast$-algebras).
\end{lemma}
\noindent {\bf Proof.} It follows  from the proof of Lemma 2 in
\cite{GZ1}.

\bigskip

Notice that $J_{00}$ and $J_{11}$ are stable under the involution
whereas $J_{01}^\ast=J_{10}$.

\bigskip

\begin{lemma} \label{identitytrasposta}
Let $\bar{A}$ be a $\ast$-simple finite dimensional algebra. Let $M
=\textrm{dim}_F \bar{A}^+$ and $L=\textrm{dim}_F \bar{A}^-$.
 Then  $\bar{A}$ does not satisfy
 $Cap^\ast_{M}[Y,X]$ and $Cap^\ast_{L}[Z,X]$.
\end{lemma}
\noindent {\bf Proof.} The result follows immediately from
\cite[Lemma 3.1]{DvS}.

\bigskip

From now on we assume that $A=M_{k}(F)+ J$, where $J=J(A)$ is the
Jacobson radical of the finite dimensional $\ast$-algebra $A$ and
$(M_k(F),t)$ is the $\ast$-algebra of matrices with transpose
involution.

\bigskip

\begin{lemma} \label{j10trasposta}
Let $M=k(k+1)/2$ and $L=k(k-1)/2$ with $k \in \mathbb{N}$, $k>0$. If
 $\Gamma^\ast_{M+1,L+1}\subseteq Id^{\ast}(A)$, then $J_{10}=J_{01}=(0)$.
\end{lemma}
\noindent {\bf Proof.} By Lemma \ref{identitytrasposta},
$M_{k}(F)$ does not satisfy the $\ast$-Capelli polynomial
$Cap^\ast_{M}[Y,X]$. Then, there exist elements $a_1^+, \ldots ,
a_M^+ \in M_{k}(F)^{+}$ and $b_1, \ldots , b_{M-1} \in  M_{k}(F)$
such that
$$
Cap^\ast_{M}(a_1^{+}, \ldots , a_M^{+};b_1, \ldots ,
b_{M-1})=e_{1,k},
$$
where the $e_{i,j}$'s are the usual matrix units. Let $d \in
J_{01},$ then $d^\ast \in J_{10}$ and $d+d^\ast \in (J_{01} \oplus
J_{10})^{+}$. Since $\Gamma^\ast_{M+1, L+1}\subseteq Id^{\ast}(A)$
we have
$$
0=Cap^\ast_{M+1}(a_1^{+}, \ldots , a_M^{+},d+d^\ast;b_1, \ldots ,
b_{M-1},e_{k,k})=de_{1,k}\pm  e_{1,k}d^\ast.
$$
Hence $de_{1,k}\pm  e_{1,k}d^\ast=0$ and, so, $de_{1,k} =\mp
e_{1,k}d^\ast \in J_{01}\cap J_{10}=(0)$. Then $d=0$, for all $d \in
J_{01}$. Thus $J_{01}=(0)$ and $J_{10}=(0)$.

\bigskip

\begin{lemma} \label{Ntrasposta}
Let $M=k(k+1)/2$ and $L=k(k-1)/2$ with $k \in \mathbb{N}$, $k>0$.
Let $J_{11} \cong M_{k}(F) \otimes_F N$, as in Lemma
\ref{radicaltrasposta}.
 If
$\Gamma_{M+1,L+1}\subseteq Id^{\ast}(A)$, then $N$ is commutative.
\end{lemma}
\noindent {\bf Proof.} Let $N$ be the finite dimensional nilpotent
$\ast$-algebra of the Lemma \ref{radicaltrasposta}. Write
$N=N^+\oplus N^-$,  where $N^+$ and $N^-$ denote the subspaces of
symmetric and skew elements of $N$ respectively. Let  $e_1^+,\ldots
, e_M^+$ be an ordered basis of $M_{k}(F)^{+}$ consisting of
symmetric matrices $e_h^+ \in \{e_{i,i}\, | \, i=1,\ldots k\}\cup \{
e_{i,j}+ e_{j,i}\, | \,  i < j , \, i,j= 1,\ldots k\} $ such that
$e_1^+= e_{1,1}$ and let $a_0,a_1, \ldots , a_M \in M_{k}(F)$ be
such that

$$
a_0e_1^+a_1\cdots a_{M-1} e_M^+ a_M=e_{1,1}
$$
and
$$
a_0e_{\sigma (1)}^+a_1\cdots a_{M-1} e_{\sigma (M)}^+ a_M=0
$$
for any $\sigma \in S_M$, $\sigma \neq id$.

Consider
$d^{+}_1, d^{+}_2 \in N^{+}$ and set $c^{+}_1=
e_{1,1} d^{+}_1$ and $c^{+}_2= e_{1,1} d^{+}_2$.
Notice that, since $N$ commutes with $M_{k}(F)$,  $c^{+}_1, c^{+}_2\in A^{+}$ and
$$
Cap^\ast_{M+2}(c^{+}_1,e_1^+,\ldots ,e_M^+, c^{+}_2; a_0, \ldots , a_M )=
$$
$$
c^{+}_1e_{1,1}c^{+}_2-c^{+}_2e_{1,1}c^{+}_1-e_{1,1}c^{+}_1c^{+}_2+$$
$$
c^{+}_2c^{+}_1e_{1,1} +e_{1,1}c^{+}_2c^{+}_1-c^{+}_1c^{+}_2e_{1,1}=
$$
$$
[c^{+}_2,c^{+}_1]e_{1,1}=[d^{+}_2,
d^{+}_1]e_{1,1}.
$$
Since $Cap^\ast_{M+1}[Y;X]\subseteq Id^{\ast}(A)$ we have $[d^{+}_1,
d^{+}_2]=0$. Thus $d^{+}_1d^{+}_2=d^{+}_2d^{+}_1$, for all $d^{+}_1,
d^{+}_2 \in N^{+}$.

Now, let $e_1^-,\ldots , e_L^-$ be an ordered basis of
$M_{k}(F)^{-}$ consisting of skew matrices $e_h^- \in  \{ e_{i,j}-
e_{j,i}\, | \,  i < j , \, i,j= 1,\ldots k\} $ such that $e_1^-=
e_{1,2}-e_{2,1}$. We consider $b_0, \ldots , b_L \in M_{k}(F)$ such
that $b_0=e_{1,1}$, $b_L=e_{k,k},$
$$
b_0e_1^-b_1 \cdots b_{L-1}e_L^-b_L=e_{1,k}
$$
and
$$
b_0e_{\tau (1)}^-b_1\cdots b_{L-1} e_{\tau (L)}^- b_L=0
$$
for all $\tau \in S_L$, $\tau \neq id$.

 Let $d^{-}_1, d^{-}_2 \in
N^{-}$ and put $c^{-}_1=  (e_{1,2}+e_{2,1})d_1^-$ and $c^{-}_2=
(e_{1,2}+e_{2,1})d_2^-$. Since $N$ commutes with $M_{k}(F)$ then
$c^{-}_1, c^{-}_2\in A^{-}$. As above we compute
$$
Cap^\ast_{L+2}(c^{-}_1,e_1^-,\ldots ,e_L^-, c^{-}_2;b_0, \ldots , b_L )=
$$
$$
[c^{-}_1,  c^{-}_2]e_{1,k}=(e_{1,1}+e_{2,2})e_{1,k}[d^{-}_1,
d^{-}_2]=e_{1,k}[d^{-}_1, d^{-}_2].
$$
Since $Cap^\ast_{L+1}[Z;X]\subseteq Id^{\ast}(A)$ we get that
$[d^{-}_1,d^{-}_2]=0$, then $d^{-}_1d^{-}_2=d^{-}_2d^{-}_1$, for all
$d^{-}_1, d^{-}_2 \in N^{-}$.

Next we show that $N^{+}$ commutes with $N^{-}$. Take $e_1^+,\ldots
, e_M^+$  an ordered basis of $M_{k}(F)^{+}$ such that $e_1^+= e_{1,1}$ and let $a_1, \ldots , a_{M} \in
M_{k}(F)$ be such that

$$
a_1e_1^+a_2\cdots a_{M} e_{M}^+=e_{1,k}
$$
and
$$
a_1e_{\rho (1)}^+a_2\cdots a_{M} e_{\rho(M)}^+ =0
$$
for any $\rho \in S_M$, $\rho \neq id$. Notice that $a_1=
a_2=e_{1,1}$. Let $d^{+}_1 \in  N^{+}$ and  $d^{-}_2 \in N^{-}$. We
set $c^{+}_1=(e_{1,2}+ e_{2,1}) d^{+}_1$ and $\bar{a}_2= e_{1,1}
d_2^{-}$. Notice that $c^{+}_1 \in A^{+}$. Then, since
$Cap^\ast_{M+1}[Y;X]\subseteq Id^{\ast}(A)$, we obtained
$$
0=Cap^\ast_{M+1}(c^{+}_1,e_1^+,\ldots ,e_M^+; a_1, \bar{a}_2, a_3, \ldots
, a_{M})=[d^{+}_1, d^{-}_2]e_{2,k}.
$$
Thus $d^{+}_1d^{-}_2=d^{-}_2d^{+}_1$, for all $d^{+}_1 \in  N^{+}$,
$d^{-}_2 \in N^{-}$ and we are done.

\bigskip

Now we are able to prove the main result about
$Id^{\ast}((M_{k}(F),t))$ and the T-$\ast$-ideal generated by the
$\ast$-Capelli polynomials $Cap^+_{\frac{k(k+1)}{2}+1}$,
$Cap^-_{\frac{k(k-1)}{2}+1}.$

\bigskip
First we prove the following
\bigskip

\begin{lemma}\label{esponentetrasposta}
Let  $M=k(k+1)/2$ and $L=k(k-1)/2$ with $k \in \mathbb{N}$, $k>0$.
Then
$$
\mathrm{exp}^\ast(\mathcal{U}^{\ast}_{M+1,L+1})=M+L=k^2=\mathrm{exp}^\ast((M_{k}(F),t)).
$$
\end{lemma}
\noindent {\bf Proof.} The exponent of
$\mathcal{U}^{\ast}_{M+1,L+1}$ is equal to the exponent of some
minimal variety lying in $\mathcal{U}^{\ast}_{M+1,L+1}$ ( for the
definition of minimal variety see \cite{GZ}).
 By \cite[Theorem 1.2]{DvS1} and Lemma \ref{capelliidentita} we have that
$$
\mathrm{exp}^\ast(\mathcal{U}^{\ast}_{M+1,L+1})=$$
$$
\mathrm{max}\{\mathrm{exp}^{\ast}(UT^\ast(A_1, \dots ,A_n)) \, | \,
UT^\ast(A_1, \dots ,A_n) \, \mathrm{satisfies} \, Cap^+_{M+1} \,
\mathrm{and}\, Cap^-_{L+1}\}.
$$

Let $d^+:=\sum_{i=1}^n \mathrm{dim}_FA_i^+$ and $d^-:=\sum_{i=1}^n
\mathrm{dim}_FA_i^-, $ then
$$
\mathrm{exp}^\ast(\mathcal{U}^{\ast}_{M+1,L+1})=
\mathrm{max}\{\mathrm{exp}^{\ast}(UT^\ast(A_1, \dots ,A_n)) \, | \,
d^++n\leq M+1 \, \mathrm{and}\, d^-+n\leq L+1\}\geq
$$
$$
\mathrm{exp}^\ast(UT^\ast (M_k(F)))=k^2=M+L.
$$

\noindent Since by Lemma \ref{esponente},
$\mathrm{exp}^\ast(\mathcal{U}^{\ast}_{M+1,L+1})\leq M+L$ then the
proof is completed.

\bigskip

\begin{theorem} \label{teorematrasposta}
Let  $M=k(k+1)/2$ and $L=k(k-1)/2$ with $k \in \mathbb{N}$, $k>0.$
Then

$$\mathcal{U}^{\ast}_{M+1,L+1}=\mathrm{var}^\ast(\Gamma^{\ast}_{M+1,L+1})=\mathrm{var}^\ast(M_{k}(F)\oplus
D'),$$

\noindent where $D'$ is a finite dimensional $\ast$-algebra such
that $\mathrm{exp}^{\ast}(D')<M+L$. In particular
$$
c^{\ast}_n(\Gamma^{\ast}_{M+1,L+1})\simeq c^{\ast}_n(M_{k}(F)).
$$
\end{theorem}

\noindent {\bf Proof.} By Lemma \ref{esponentetrasposta} we have
that $\mathrm{exp}^{\ast}(\mathcal{U}^{\ast}_{M+1,L+1})=M+L$.

Let $A=A^+\oplus A^-$ be a generating $\ast$-algebra of
$\mathcal{U}^{\ast}_{M+1,L+1}$. As remarked before
we can assume that $A$ is  finite dimensional. Thus, by Corollary
\ref{decomposizionefinita}, there exists a finite
 number of reduced $\ast$-algebras $B_1,\ldots , B_s$ and a finite dimensional $\ast$-algebra $D'$ such that

$$
\mathcal{U}^{\ast}_{M+1,L+1}=\mathrm{var}^\ast(A)=\mathrm{var}^\ast(B_1\oplus
\cdots \oplus B_s \oplus D'). \eqno (1)
$$

\noindent Moreover $$\mathrm{exp}^{\ast}(B_1)=\cdots
=\textrm{exp}^{\ast}(B_s)=
\textrm{exp}^{\ast}(\mathcal{U}^{\ast}_{M+1,L+1})=M+L$$
 and \,

 $$\textrm{exp}^{\ast}(D')< \textrm{exp}^{\ast}(\mathcal{U}^{\ast}_{M+1,L+1})= M+L.$$

 Next, we
analyze the structure of a finite dimensional reduced $\ast$-algebra
 $R$  such that
$\textrm{exp}^{\ast}(R)= M+L=
\textrm{exp}^{\ast}(\mathcal{U}^{\ast}_{M+1,L+1})$ and
$\Gamma^\ast_{M+1,L+1}\subseteq Id^{\ast}(R) $. We can write

 $$R=R_1
\oplus \cdots \oplus R_q + J ,$$ where $R_i$ are simple
$\ast$-subalgebras of $R$, $J=J(R)$ is the Jacobson radical of $R$
and $R_1J\cdots J R_q\neq 0$.

Recall that every $\ast$-algebra $R_i$  is isomorphic to one of the
following algebras :$(M_{k_i}(F),t)$ or  $(M_{2m_i}(F),s)$  or
$(M_{h_i}(F)\oplus M_{h_i}(F)^{op},exc).$

Let $t_1$ be the number of $\ast$-algebras $R_i$ of the first type,
 $t_2$ the number of $\ast$-algebras $R_i$ of the second type
and  $t_3$  the number of $R_i$ of the third type, with
$t_1+t_2+t_3=q$.

By \cite[Theorem 4.5]{DvLs} and \cite[Proposition 4.7]{DvLs} there exists a $\ast$-algebra $\overline{R}$ isomorphic to the
$\ast$-algebra $UT^\ast(R_1, \ldots, R_q)$ such that

$$\mathrm{exp}^{\ast}(R)=\mathrm{exp}^{\ast}(\overline{R})=\mathrm{exp}^{\ast}(UT^\ast(R_1,
\ldots, R_q)).$$

\noindent Let observe that
$$
k^2=M+L=\mathrm{exp}^{\ast}(R)=\mathrm{exp}^{\ast}(\overline{R})=\mathrm{exp}^{\ast}(UT^\ast(R_1,
\ldots, R_q)) =
$$
$$ \mathrm{dim}_FR_1+ \cdots + \mathrm{dim}_F R_q=
k_1^2 + \cdots + k_{t_1}^2+ (2 m_1)^2+ \cdots + (2 m_{t_2})^2 +
2h_1^2 + \cdots + 2h_{t_3}^2.
$$

\noindent Let $d^{\pm}=\mathrm{dim}_F(R_1\oplus \cdots \oplus
R_q)^{\pm}$ then
$$
d^++d^-=d=\mathrm{dim}_F(R_1\oplus \cdots \oplus
R_q)=\mathrm{exp}^{\ast}(\overline{R})=M+L.
$$
By \cite[Lemma 3.2]{DvS1} $\overline{R}$ does not satisfy the $\ast$-Capelli polynomials
$Cap^\ast_{d^++q-1}[Y;X]$ and $Cap^\ast_{d^-+q-1}[Z;X]$, but
$\overline{R}$ satisfies $Cap^\ast_{M+1}[Y;X]$ and
$Cap^\ast_{L+1}[Z;X]$. Thus $d^++q-1\leq M$ and  $d^-+q-1\leq L$.
Hence $d^++ d^- + 2q-2\leq M+L$. Since $d^++d^-= M+L$ we obtain that
$2q-2= 0$ and so $1=q= t_1+t_2 +t_3$. Since $t_1$, $t_2$ and $t_3$
are nonnegative integers, we have the following three possibilities
\begin{enumerate}
  \item $t_1=1$ and $t_2=t_3=0$;
  \item $t_2=1$ and $t_1=t_3=0$
  \item $t_3=1$ and $t_1=t_2=0$.
\end{enumerate}

\noindent If $t_2=1$, then $R=(M_{2m}(F), s)+J$ and
$\mathrm{exp}^\ast(R)=4m^2$. Thus
$$k^2=M+L=\mathrm{exp}^{\ast}(R)=4m^2$$ and so $k=2m$. By hypothesis $R$
satisfies  $Cap^\ast_{L+1}[Z;X]$ but, since $Id^\ast(R)\subseteq
Id^\ast(\overline{R})$, $R$  does not satisfy $Cap^\ast_{d^-}[Z;X],$
where $d^- = m(2m+1).$ It follows that
$$L+1=k(k-1)/2+1=m(2m-1)+1<m(2m+1)=d^-,$$ for $m\geq 1,$ and this is a
contradiction.

\noindent Let assume  $t_3=1$. Then $R=(M_{h}(F)\oplus
M_{h}(F)^{op})+J$ and $\mathrm{exp}^{\ast}(R)=2h^2$. Thus
$k^2=M+L=\mathrm{exp}^\ast(R)=2h^2$ and this is impossible.

\noindent Then $t_1 =1$ and in this case
$$
R\cong M_k(F)+J.
$$
From Lemmas \ref{radicaltrasposta}, \ref{j10trasposta},
\ref{Ntrasposta} we obtain

$$
R\cong (M_k(F)+J_{11})\oplus J_{00}\cong (M_k(F)\otimes
N^\sharp)\oplus J_{00}
$$

\noindent where $N^\sharp$ is the algebra obtained from $N$ by
adjoining a unit element. Since $N^\sharp$ is commutative, it
follows that $M_k(F) + J_{11}$ and $M_k(F)$ satisfy the same
$\ast$-identities. Thus var$^\ast$(R)=var$^\ast(M_k(F)\oplus J_{00})$ with
$J_{00}$ a finite dimensional nilpotent $\ast$-algebra. Hence,
recalling the decomposition given in (1), we get

$$\mathcal{U}^{\ast}_{M+1,L+1}=\mathrm{var}^{\ast}(\Gamma_{M+1,L+1})=\mathrm{var}^{\ast}(M_k(F)\oplus
D'),$$

\noindent where $D'$ is a finite dimensional $\ast$-algebra with
exp$^{\ast}(D')<M+L$. Then, from Corollary \ref{codimensioni} we
have

$$
c^{\ast}_n(\Gamma_{M+1,L+1})\simeq c^{\ast}_n(M_k(F))
$$

\noindent and the theorem is proved.

\bigskip

\section{Asymptotics for $\mathcal{U}^\ast_{m(2m-1)+1,m(2m+1)+1}$ and $(M_{2m}(F),s)$}

\bigskip

Throughout this section we assume that $A=M_{2m}(F)+ J$, where
$J=J(A)$ is the Jacobson radical of the finite dimensional
$\ast$-algebra $A$ and $(M_{2m}(F),s)$ is the algebra of matrices
with symplectic involution.

\begin{lemma} \label{j10simplettica}
Let $M=m(2m-1)$ and $L=m(2m+1)$ with $m \in \mathbb{N}$, $m>0$. If
 $\Gamma^\ast_{M+1,L+1}\subseteq Id^{\ast}(A)$, then $J_{10}=J_{01}=(0)$.
\end{lemma}
\noindent {\bf Proof.} By Lemma \ref{identitytrasposta},
$(M_{2m}(F),s)$ does not satisfy the $\ast$-Capelli polynomial
$Cap^\ast_{M}[Y,X]$. Also there exist elements $a_1^+, \ldots ,
a_M^+ \in M_{2m}(F)^{+}$ and $b_1, \ldots , b_{M-1} \in M_{2m}(F)$
such that
$$
Cap^\ast_{M}(a_1^{+}, \ldots , a_M^{+};b_1, \ldots ,
b_{M-1})=e_{1,2m},
$$
where the $e_{i,j}$'s are the usual matrix units.  Let $d \in
J_{01},$ then $d+d^\ast \in (J_{01} \oplus J_{10})^{+}$. Since
$\Gamma^\ast_{M+1, L+1}\subseteq Id^{\ast}(A)$ we have
$$
0=Cap^\ast_{M+1}(a_1^{+}, \ldots , a_M^{+},d+d^\ast;b_1, \ldots ,
b_{M-1},e_{2m,2m})=de_{1,2m}\pm  e_{1,2m}d^\ast.
$$
Hence $de_{1,2m} =\mp e_{1,2m}d^\ast \in J_{01}\cap J_{10}=(0)$ and
so $d=0$, for all $d \in J_{01}$. Thus $J_{01}=(0)=J_{10}$ and we
are done.

\bigskip

\begin{lemma} \label{Nsimplettica}
Let $M=m(2m-1)$ and $L=m(2m+1)$ with $m \in \mathbb{N}$, $m>0$. Let
$J_{11} \cong M_{2m}(F) \otimes_F N$, as in Lemma
\ref{radicaltrasposta}.  If $\Gamma_{M+1,L+1}\subseteq
Id^{\ast}(A)$, then $N$ is commutative.
\end{lemma}
\noindent {\bf Proof.} Let $N$ be the finite dimensional nilpotent
$\ast$-algebra of the Lemma \ref{radicaltrasposta}.  As in Lemma
\ref{Ntrasposta} we can consider an ordered basis of $M_{2m}(F)^{+}$
consisting of symmetric matrices   $e_h^+ \in \{e_{i,j}+e_{m+j,m+i}\,
| \, i,j=1,\ldots m\}\cup \{ e_{i,m+j}- e_{j,m+i}\, | \,  1 \leq i <
j \leq m\} \cup \{ e_{m+i,j}- e_{m+j,i}\, | \,  1 \leq i < j \leq
m\}$ with $e_1^+= e_{1,1}+e_{m+1,m+1}$ and  $a_0,a_1, \ldots , a_M
\in M_{2m}(F)$ such that
$$
a_0e_1^+a_1\cdots a_{M-1} e_M^+ a_M=e_{1,1}
$$
and
$$
a_0e_{\sigma (1)}^+a_1\cdots a_{M-1} e_{\sigma (M)}^+ a_M=0
$$

\noindent for any $\sigma \in S_M$, $\sigma \neq id$. Consider
$d^{+}_1, d^{+}_2 \in N^{+}$ and set $c^{+}_1= (e_{1,1}+e_{m+1,m+1})
d^{+}_1$ and $c^{+}_2= (e_{1,1}+e_{m+1,m+1}) d^{+}_2$. Since $N$
commutes with $M_{2m}(F)$,  $c^{+}_1, c^{+}_2\in A^{+}$ and we
obtain
$$
Cap^\ast_{M+2}(c^{+}_1,e_1^+,\ldots ,e_M^+, c^{+}_2; a_0, \ldots , a_M )=
[c^{+}_2,c^{+}_1]e_{1,1}=[d^{+}_2, d^{+}_1]e_{1,1}.
$$

\noindent Since $Cap^\ast_{M+1}[Y;X]\subseteq Id^{\ast}(A)$ we have
$d^{+}_1d^{+}_2=d^{+}_2d^{+}_1$, for all $d^{+}_1, d^{+}_2 \in
N^{+}$.

Now, let $e_1^-,\ldots , e_L^-$ be an ordered basis of $M_{2m}(F)^{-}$
consisting of skew matrices, $e_h^- \in \{e_{i,j}-e_{m+j,m+i}\, | \,
i,j=1,\ldots m\}\cup \{ e_{i,m+j}+ e_{j,m+i}\, | \,  1 \leq i < j
\leq m\}\cup \{ e_{i,m+i}\, | \, i=1,\ldots , m\} \cup \{ e_{m+i,j}+
e_{m+j,i}\, | \,  1 \leq i < j \leq m\}\cup \{ e_{m+i,i}\, | \,
i=1,\ldots , m\}$ such that $e_1^-= e_{1,1}-e_{m+1,m+1}$ .

We consider $b_0, \ldots , b_L \in M_{2m}(F)$ such that
$b_0=e_{1,1}$, $b_L=e_{2m,2m}$
$$
b_0e_1^-b_1 \cdots b_{L-1}e_L^-b_L=e_{1,2m}
$$
and
$$
b_0e_{\tau (1)}^-b_1\cdots b_{L-1} e_{\tau (L)}^- b_L=0
$$
for all $\tau \in S_L$, $\tau \neq id$. Let $d^{-}_1, d^{-}_2 \in
N^{-}$. Let $c^{-}_1= (e_{1,1}+e_{m+1,m+1}) d^{-}_1$ and $c^{-}_2=
(e_{1,1}+e_{m+1,m+1}) d^{-}_2$. Since $N$ commutes with $M_{2m}(F)$
then $c^{-}_1, c^{-}_2\in A^{-}$. As above we obtain

$$
Cap^\ast_{L+2}(c^{-}_1,e_1^-,\ldots ,e_L^-, c^{-}_2;b_0, \ldots , b_L )=
$$
$$
[c^{-}_1,  c^{-}_2]e_{1,2m}=(e_{1,1}+e_{m+1,m+1})e_{1,2m}[d^{-}_1,
d^{-}_2]=e_{1,2m}[d^{-}_1, d^{-}_2].
$$

\noindent Since $Cap^\ast_{L+1}[Z;X]\subseteq Id^{\ast}(A)$ we get
that $[d^{-}_1,d^{-}_2]=0$ and so $d^{-}_1d^{-}_2=d^{-}_2d^{-}_1$,
for all $d^{-}_1, d^{-}_2 \in N^{-}$.

Next we show that $N^{+}$ commutes with $N^{-}$.

 Take $e_1^+,\ldots
, e_M^+$  an ordered basis of $M_{2m}(F)^{+}$  such that $e_1^+= e_{1,1}+e_{m+1,m+1}$. Let $a_1, \ldots ,
a_{M} \in M_{2m}(F)$ be such that
$$
a_1e_1^+a_2\cdots a_{M} e_{M}^+=e_{1,2m}
$$
and
$$
a_1e_{\rho (1)}^+a_2\cdots a_{M} e_{\rho(M)}^+ =0
$$
for any $\rho \in S_M$, $\rho \neq id$. Notice that $a_1=
a_2=e_{1,1}$. Let $d^{+}_1 \in  N^{+}$ and  $d^{-}_2 \in N^{-}$. We
set $c^{+}_1=(e_{1,2}+ e_{m+2,m+1}) d^{+}_1$ and $\bar{a}_2= e_{1,1}
d_2^{-}$, then  $c^{+}_1 \in A^{+}$. Since
$Cap^\ast_{M+1}[Y;X]\subseteq Id^{\ast}(A)$, we obtain

$$
0=Cap^\ast_{M+1}(c^{+}_1,e_1^+,\ldots ,e_M^+; a_1, \bar{a}_2, a_3, \ldots
, a_{M})=[d^{+}_1, d^{-}_2]e_{2,2m}.
$$

\noindent  Thus $d^{+}_1d^{-}_2=d^{-}_2d^{+}_1$, for all $d^{+}_1
\in N^{+}$, $d^{-}_2 \in N^{-}$ and we are done.

\bigskip

\begin{lemma}\label{esponentesimplettica}
Let  $M=m(2m-1)$ and $L=m(2m+1)$ with $m \in \mathbb{N}$, $m>0$.
Then

$$
\mathrm{exp}^\ast(\mathcal{U}^{\ast}_{M+1,L+1})=M+L=4m^2=\mathrm{exp}^\ast((M_{2m}(F),s)).
$$
\end{lemma}
\noindent {\bf Proof.} The proof is the same  of that of Lemma
\ref{esponentetrasposta}.

\bigskip

Now we are able to prove the following

\bigskip

\begin{theorem} \label{teoremasimplettica}
Let $M=m(2m-1)$ and $L=m(2m+1)$ with $m \in \mathbb{N}$, $m>0$.
Then

$$\mathcal{U}^{\ast}_{M+1,L+1}=\mathrm{var}^\ast(\Gamma_{M+1,L+1})=\mathrm{var}^\ast(M_{2m}(F)\oplus
D''),$$

\noindent  where $D''$ is a finite dimensional $\ast$-algebra such
that $\mathrm{exp}^{\ast}(D'')<M+L$. In particular

$$
c^{\ast}_n(\Gamma_{M+1,L+1})\simeq c^{\ast}_n(M_{2m}(F)).
$$
\end{theorem}

\noindent {\bf Proof.}  The first part of the proof follows step by step
 that of Theorem \ref{teorematrasposta} and we obtain

$$
\mathcal{U}^{\ast}_{M+1,L+1}=\mathrm{var}^\ast(B_1\oplus
\cdots \oplus B_s \oplus D''),
$$

\noindent where  $B_1,\ldots , B_s$ are reduced $\ast$-algebras and
$D''$ is a finite dimensional
 $\ast$-algebra such that $\mathrm{exp}^\ast(D'')<\mathrm{exp}^\ast(B_i)=M+L$, for all $i=1,\ldots s$.

Let $R$ be a finite
dimensional reduced $\ast$-algebra  with $\mathrm{exp}^\ast(R)=M+L=\mathrm{exp}^\ast(\mathcal{U}_{M+1,L+1})$ and
$\Gamma^\ast_{M+1,L+1}\subseteq Id^{\ast}(R)$.  We can write $R=R_1
\oplus \cdots \oplus R_q + J$, where $R_i$ are simple
$\ast$-subalgebras of $R$, $J=J(R)$ is the Jacobson radical of $R$.
Let $t_1$ be the number of $\ast$-algebras $R_i$ isomorphic to $(M_{k_i}(F),t)$,
 $t_2$  the number of $\ast$-algebras $R_i$ isomorphic to $(M_{2m_i}(F),s)$
and let $t_3$ be the number of $R_i$ isomorphic to
$(M_{h_i}(F)\oplus M_{h_i}(F)^{op},exc),$  where $t_1+t_2+t_3=q$.
Hence, as in Theorem \ref{teorematrasposta}, there exists a
$\ast$-algebra $\overline{R}$ isomorphic to the $\ast$-algebra
$UT^\ast(R_1, \ldots, R_q)$ such that
$\mathrm{exp}^{\ast}(R)=\mathrm{exp}^{\ast}(\overline{R})=\mathrm{exp}^{\ast}(UT^\ast(R_1,
\ldots, R_q))$ and
$$
4m^2=M+L=
\mathrm{exp}^{\ast}(R)=\mathrm{exp}^{\ast}(\overline{R})=\mathrm{exp}^{\ast}(UT^\ast(R_1,
\ldots, R_q))=
$$
$$
k_1^2 + \cdots + k_{t_1}^2+ (2 m_1)^2+ \cdots + (2 m_{t_2})^2 +
2h_1^2 + \cdots + 2h_{t_3}^2.
$$

As in the proof of the Theorem \ref{teorematrasposta} we have $q=1$ and so we obtain only three possibilities: $t_1=1$ and $t_2=t_3=0$ or $t_2=1$ and $t_1=t_3=0$ or $t_3=1$ and $t_1=t_2=0$.

If $t_1=1$, then $R=(M_{k}(F), t)+J$ and
$\mathrm{exp}^\ast(R)=k^2$. Thus $4m^2=M+L=\mathrm{exp}^{\ast}(R)=k^2$ and so $k=2m$. Hence $R$
satisfies  $Cap^\ast_{M+1}[Y;X]$ and, since $Id^\ast(R)\subseteq Id^\ast(\overline R)$, $R$  does not satisfy $Cap^\ast_{d^+}[Y;X],$
where $d^+ = k(k+1)/2$. It follows that
$M+1=m(2m-1)+1<m(2m+1)=k(k+1)/2=d^+$ for $m\geq 1,$ and this is a
contradiction.

Let assume  $t_3=1$. Then $R=(M_{h}(F)\oplus M_{h}(F)^{op})+J$ and
$\mathrm{exp}^{\ast}(R)=2h^2$. Thus  $4m^2=M+L=\mathrm{exp}^\ast(R)=2h^2$ and
this is impossible.

Finally, let $t_2 =1$ and $t_1=t_3=0$. Then $R\cong M_{2m}(F)+J$. By
Lemmas \ref{radicaltrasposta}, \ref{j10simplettica},
\ref{Nsimplettica} we obtain

$$
R\cong (M_{2m}(F)+J_{11})\oplus J_{00}\cong (M_{2m}(F)\otimes
N^\sharp)\oplus J_{00}
$$

\noindent where $N^\sharp$ is the algebra obtained from $N$ by
adjoining a unit element. Since $N^\sharp$ is commutative, it
follows that $M_{2m}(F) + J_{11}$ and $M_{2m}(F)$ satisfy the same
$\ast$-identities. Thus var$^\ast$(R)=var$^\ast(M_{2m}(F)\oplus
J_{00})$ with $J_{00}$ a finite dimensional nilpotent
$\ast$-algebra. Hence by (1) we get

$$\mathcal{U}^{\ast}_{M+1,L+1}=\mathrm{var}^{\ast}(\Gamma_{M+1,L+1})=\mathrm{var}^{\ast}(M_{2m}(F)\oplus
D''),$$

\noindent where $D''$ is a finite dimensional $\ast$-algebra with
exp$^{\ast}(D'')<M+L$. Then, from Corollary \ref{codimensioni}, we
have
$$
c^{\ast}_n(\Gamma_{M+1,L+1})\simeq c^{\ast}_n(M_{2m}(F))
$$

\noindent and the theorem is proved.

\bigskip

\section{Asymptotics for $\mathcal{U}^\ast_{h^2+1,h^2+1}$ and $(M_{h}(F)\oplus M_{h}(F)^{op}, exc)$}

\bigskip
Throughout this section we assume that $A=(M_{h}(F)\oplus M_{h}(F)^{op})+ J$, where
$J=J(A)$ is the Jacobson radical of the finite dimensional
$\ast$-algebra $A$ and $(M_{h}(F)\oplus M_{h}(F)^{op}, exc)$ is the  direct
sum of the algebra of $h \times h$ matrices and the opposite
algebra with the exchange involution.

We start with the following lemmas

\begin{lemma} \label{j10scambio}
Let $M=L=h^2$ with $h \in \mathbb{N}$, $h>0$. If
 $\Gamma^\ast_{M+1,L+1}\subseteq Id^{\ast}(A)$, then $J_{10}=J_{01}=(0)$.
\end{lemma}
\noindent {\bf Proof.} By Lemma \ref{identitytrasposta},
$(M_{h}(F)\oplus M_{h}(F)^{op}, exc)$ does not satisfy the $\ast$-Capelli polynomial
$Cap^\ast_{M}[Y,X]$. Also, there exist elements $a_1^+, \ldots ,
a_M^+ \in (M_{h}(F)\oplus M_{h}(F)^{op})^{+}$ and $b_1, \ldots , b_{M-1} \in M_{h}(F)\oplus M_{h}(F)^{op}$
such that
$$
Cap^\ast_{M}(a_1^{+}, \ldots , a_M^{+};b_1, \ldots ,
b_{M-1})=\tilde{e}_{1,h},
$$
where the $\tilde{e}_{i,j}=(e_{i,j}, e_{j,i})$ and $e_{i,j}$'s are
the usual matrix units.  Let $d \in J_{01},$ then $d+d^\ast \in
(J_{01} \oplus J_{10})^{+}$. Since $\Gamma^\ast_{M+1, L+1}\subseteq
Id^{\ast}(A)$ it follows
$$
0=Cap^\ast_{M+1}(a_1^{+}, \ldots , a_M^{+},d+d^\ast;b_1, \ldots ,
b_{M-1},\tilde{e}_{h,h})\tilde{e}_{h,h}=d\tilde{e}_{1,h}.
$$
Hence $d=0$, for all $d \in J_{01}$. Thus $J_{01}=(0)$ and
$J_{10}=(0)$.

\bigskip

\begin{lemma} \label{Nscambio}
Let $M=L=h^2$  with $h \in \mathbb{N}$, $h>0$. Let
$J_{11} \cong (M_{h}(F)\oplus M_{h}(F)^{op}) \otimes_F N$, as in Lemma
\ref{radicaltrasposta}.  If $\Gamma_{M+1,L+1}\subseteq
Id^{\ast}(A)$, then $N$ is commutative.
\end{lemma}
\noindent {\bf Proof.} Let $N$ be the finite dimensional nilpotent
$\ast$-algebra of the Lemma \ref{radicaltrasposta}.  Let now $v_1^+, \ldots , v_M^+$ an ordered basis of $(M_{h}(F)\oplus M_{h}(F)^{op})^+$ consisting of all $e_{i,j}^+=(e_{i,j},e_{i,j})$ such that $v_1^+=e_{1,1}^+$ and let $a_0, \dots , a_M \in M_{h}(F)\oplus M_{h}(F)^{op}$ be such that
$$
a_0 v_1^+ a_1 \cdots a_{M-1}v_M^+a_M=e_{1,1}^+
$$
and
$$
a_0 v_{\sigma(1)}^+ a_1 \cdots a_{M-1}v_{\sigma(M)}^+a_M=(0,0)
$$
for any $\sigma \in S_M$, $\sigma \neq id$.

Now let $d_1^+$, $d_2^+\in N^+$ and set $c_1^+=e_{1,1}^+d_1^+$ and
$c_2^+=e_{1,1}^+d_2^+$. Recalling that $N$ commutes with
$M_{h}(F)\oplus M_{h}(F)^{op}$ we have that $c_1^+$, $c_2^+\in A^+$
and, as in Lemma \ref{Ntrasposta},
$$
Cap^\ast_{M+2}(c^{+}_1,v_1^+,\ldots ,v_M^+, c^{+}_2; a_0, \ldots , a_M )=
[c^{+}_2,c^{+}_1]e_{1,1}^+=[d^{+}_2,d^{+}_1]e_{1,1}^+.
$$

\noindent Since $Cap^\ast_{M+1}[Y;X]\subseteq Id^\ast(A)$ it follows
$d_1^+d_2^+=d_2^+d_1^+$, for all $d_1^+$, $d_2^+\in N^+$.

Now, let $v_1^-, \ldots , v_L^-$ be an ordered basis of
$(M_{h}(F)\oplus M_{h}(F)^{op})^-$ consisting of all
$e_{i,j}^-=(e_{i,j},-e_{i,j})$ such that $v_1^-=e_{1,1}^-$ and let
$b_0, \dots , b_L \in M_{h}(F)\oplus M_{h}(F)^{op}$ be such that
$$
b_0 v_1^+ b_1 \cdots b_{L-1}v_L^+b_L=\left\{
                                       \begin{array}{ll}
                                         e_{1,1}^+, & \hbox{if L is even;} \\
                                         e_{1,1}^-, & \hbox{if L is odd.}
                                       \end{array}
                                     \right.
$$
and
$$
b_0 v_{\sigma(1)}^+ b_1 \cdots b_{L-1}v_{\sigma(L)}^+b_L=(0,0)
$$
for any $\sigma \in S_M$, $\sigma \neq id$.
Now, we take $d_1^-$, $d_2^-\in N^-$ and set $c_1^-=e_{1,1}^+d_1^-$ and $c_2^-=e_{1,1}^+d_2^-$.
Since $N$ commutes with $M_{h}(F)\oplus M_{h}(F)^{op}$, we have that $c_1^-$, $c_2^-\in A^-$.
Thus, since $Cap^\ast_{L+1}[Z,X]\subseteq Id^\ast(A)$, we obtain
$$
0\! = \! Cap^\ast_{L+2}(c^{-}_1,v_1^-,\ldots ,v_L^-, c^{-}_2;b_0, \ldots , b_L)\! = \! \left\{
                                       \begin{array}{ll}
                                         [c^{-}_2,c^{-}_1]e_{1,1}^+= [d^{-}_2,d^{-}_1]e_{1,1}^+, & \!\!\hbox{if L is even;} \\
                                         \hbox{$[c^{-}_2,c^{-}_1]$}e_{1,1}^-=[d^{-}_2,d^{-}_1]e_{1,1}^-, &\!\! \hbox{if L is odd.}
                                       \end{array}
                                     \right.
$$
In conclusion $[d^{-}_2,d^{-}_1]=0$, for all $d_1^-$, $d_2^-\in N^-$.

Finally, we consider $d_1^+ \in N^+$, $d_2^-\in N^-$ and set $c^+_1=e_{1,1}^+d_1^+$  and $c^+_2=e_{1,1}^-d_2^-$. As above, we have that $c^+_1$, $c^+_2 \in A^+$. Then
$$
0=Cap^\ast_{M+2}(c^{+}_1,v_1^+,\ldots ,v_M^+, c^{+}_2;a_0, \ldots , a_M)=[c^{+}_2,c^{+}_1]e_{1,1}^+= [d^{-}_2,d^{+}_1]e_{1,1}^-
$$
and $[d^{-}_2,d^{+}_1]=0$ for all $d_1^+ \in N^+$, $d_2^-\in N^-$
and the lemma is proved.
\bigskip

\begin{lemma}\label{esponentescambio}
Let $M=L=h^2$  with $h \in \mathbb{N}$, $h>0$.
Then
$$
\mathrm{exp}^\ast(\mathcal{U}^{\ast}_{M+1,L+1})=M+L=2h^2=\mathrm{exp}^\ast((M_{h}(F)\oplus M_{h}(F)^{op},exc)).
$$
\end{lemma}
\noindent {\bf Proof.} The proof is the same  of that of Lemma
\ref{esponentetrasposta}.

\bigskip

\begin{theorem} \label{teoremascambio}
Let $M=L=h^2$  with $h \in \mathbb{N}$, $h>0$.
Then

$$\mathcal{U}^{\ast}_{M+1,L+1}=\mathrm{var}^\ast(\Gamma_{M+1,L+1})=\mathrm{var}^\ast((M_{h}(F)\oplus M_{h}(F)^{op})\oplus
D'''),$$

\noindent  where $D'''$ is a finite dimensional $\ast$-algebra such
that $\mathrm{exp}^{\ast}(D''')<M+L$. In particular
$$
c^{\ast}_n(\Gamma_{M+1,L+1})\simeq c^{\ast}_n(M_{h}(F)\oplus M_{h}(F)^{op}).
$$
\end{theorem}

\noindent {\bf Proof.}  The  proof proceeds as in
 that of Theorem \ref{teorematrasposta}. Let $R$ be a finite
dimensional reduced $\ast$-algebra  such that
$\mathrm{exp}^\ast(R)=M+L=\mathrm{exp}^\ast(\mathcal{U}^\ast_{M+1,L+1})$
and $\Gamma^\ast_{M+1,L+1}\subseteq Id^{\ast}(R)$.  We can write
$R=R_1 \oplus \cdots \oplus R_q + J$, where $R_i$ are simple
$\ast$-subalgebras of $R$, $J=J(R)$ is the Jacobson radical of $R$.
If $t_1$ denotes the number of algebras $R_i$ isomorphic to
$(M_{k_i}(F),t)$,
 $t_2$  the number of algebras $R_i$ isomorphic to $(M_{2m_i}(F),s)$
and  $t_3$  the number of $R_i$ isomorphic to $(M_{h_i}(F)\oplus
M_{h_i}(F)^{op},exc)$, then we have $t_1+t_2+t_3=q$. As in Theorem
\ref{teorematrasposta}, there exists a $\ast$-algebra $\overline{R}$
isomorphic to the $\ast$-algebra $UT^\ast(R_1, \ldots, R_q)$ such
that
$\mathrm{exp}^{\ast}(R)=\mathrm{exp}^{\ast}(\overline{R})=\mathrm{exp}^{\ast}(UT^\ast(R_1,
\ldots, R_q))$. Then

$$
2h^2=M+L=
\mathrm{exp}^{\ast}(R)=\mathrm{exp}^{\ast}(\overline{R})=\mathrm{exp}^{\ast}(UT^\ast(R_1,
\ldots, R_q))=
$$
$$
k_1^2 + \cdots + k_{t_1}^2+ (2 m_1)^2+ \cdots + (2 m_{t_2})^2 +
2h_1^2 + \cdots + 2h_{t_3}^2
$$

\noindent and we have only three possibilities: $t_1=1$ and
$t_2=t_3=0$ or $t_2=1$ and $t_1=t_3=0$ or $t_3=1$ and $t_1=t_2=0$.

If $t_1=1$, then $R=(M_{k}(F), t)+J$ and
$\mathrm{exp}^\ast(R)=k^2$. Thus $2h^2=M+L=\mathrm{exp}^{\ast}(R)=k^2$ and this is a contradiction.

If $t_2=1$. Then $R=(M_{2m}(F),s)+J$ and
$\mathrm{exp}^{\ast}(R)=4m^2$. Thus  $2h^2=M+L=\mathrm{exp}^\ast(R)=4m^2$ and so $h^2=2m^2$, contradiction.

Finally, let $t_3 =1$ and $t_1=t_2=0$. Then $R\cong (M_{h}(F)\oplus M_{h}(F)^{op})+J$.
 As in Theorem \ref{teorematrasposta}, by Lemmas \ref{radicaltrasposta}, \ref{j10scambio}, \ref{Nscambio} we obtain

$$
R\cong ((M_{h}(F)\oplus M_{h}(F)^{op})+J_{11})\oplus J_{00}\cong
((M_{h}(F)\oplus M_{h}(F)^{op})\otimes N^\sharp)\oplus J_{00}
$$

\noindent where $N^\sharp$ is the algebra obtained from $N$ by
adjoining a unit element. Since $N^\sharp$ is commutative, we have
that $(M_{h}(F)\oplus M_{h}(F)^{op}) + J_{11}$ and $M_{h}(F)\oplus
M_{h}(F)^{op}$ satisfy the same $\ast$-identities. Thus
var$^\ast$(R)=var$^\ast((M_{h}(F)\oplus M_{h}(F)^{op})\oplus
J_{00})$ with $J_{00}$ a finite dimensional nilpotent
$\ast$-algebra. We get

$$
\mathcal{U}^{\ast}_{M+1,L+1}=\mathrm{var}^{\ast}(\Gamma_{M+1,L+1})=\mathrm{var}^{\ast}((M_{h}(F)\oplus M_{h}(F)^{op})\oplus
D'''),$$

\noindent where $D'''$ is a finite dimensional $\ast$-algebra with
exp$^{\ast}(D''')<M+L$. Then, from Corollary \ref{codimensioni}, we
have
$$
c^{\ast}_n(\Gamma_{M+1,L+1})\simeq c^{\ast}_n((M_{h}(F)\oplus M_{h}(F)^{op}))
$$

\noindent and the proof is completed.

\bigskip

\bibliographystyle{amsplain}

\end{document}